\magnification=\magstep1
\hsize=5.9truein
\voffset=.5truein
\input amssym.def
\input amssym.tex
\def\newline{\hfill\break}
\def\scong{{\scriptstyle\|}\lower.2ex\hbox{$\wr$}}
\def\Z{{\Bbb Z}}

\def\Q{{\Bbb Q}}

\def\Jac{\mathop{\rm Jac}\nolimits}
\def\Tor{\mathop{\rm Tor}\nolimits}
\def\Div{\mathop{\rm Div}\nolimits}

\def\Sym{\mathop{\rm Sym}\nolimits}

\def\Hom{\mathop{\rm Hom}\nolimits}

\def\Pic{\mathop{\rm Pic}\nolimits}
\def\Gal{\mathop{\rm Gal}\nolimits}
\def\Br{\mathop{\rm Br}\nolimits}
\def\ram{\mathop{\rm ram}\nolimits}
\def\Spec{\mathop{\rm Spec}\nolimits}

\def\rtimes{\mathop{\times\!\!{\raise.2ex\hbox{$\scriptscriptstyle|$}}}
    \nolimits}
\def\proof{\noindent{\it Proof.}\quad}
\def\blackbox{\hbox{\vrule width6pt height7pt depth1pt}}
\outer\def\Demo #1. #2\par{\medbreak\noindent {\it#1.\enspace}
    {\rm#2}\par\ifdim\lastskip<\medskipamount\removelastskip
    \penalty55\medskip\fi}
\def\qed{~\hfill\blackbox\medskip}
\overfullrule=0pt
\def\Br{\mathop{\rm Br}\nolimits}

\def\hangbox to #1 #2{\vskip1pt\hangindent #1\noindent \hbox to #1{#2}$\!\!$}

\pageno=0
\footline{\ifnum\pageno=0\hfill\else\hss\tenrm\folio\hss\fi}
\topinsert\vskip1.8truecm\endinsert
\centerline{\bf Tensor Products of Division Algebras and Fields}
\vskip6pt
$${\vbox{\halign{\hfil\hbox{#}\hfil\qquad&\hfil\hbox{#}\hfil\cr
Louis Rowen\cr
Department of Mathematics\cr
Bar-Ilan University\cr
Ramat Gan, 52900, Israel\cr}}}$$
\vskip6pt
$${\vbox{\halign{\hfil\hbox{#}\hfil\qquad&\hfil\hbox{#}\hfil\cr
David J. Saltman\cr
Center for Communications Research\cr
Princeton, NJ 08540\cr}}}$$
\vskip16pt
\vskip3pt
{\narrower\smallskip\noindent
{\bf Abstract}
This paper began as an investigation of the question of whether
$D_1 \otimes_F D_2$ is a domain where the $D_i$ are division algebras
and $F$ is an algebraically closed field 
contained in their centers. We present an example where the answer
is ``no'', and also study the Picard group and Brauer group
properties of $F_1 \otimes_F F_2$ where the $F_i$ are fields.
Finally, as part of our example, we have results about division
algebras and Brauer groups over curves. Specifically, we give a
splitting criterion for certain Brauer group elements on the
product of two curves over $F$.

\bigskip

\noindent AMS Subject Classification:  16K20, 14F22, 16K50, 12G05, 14C22
\medskip

\noindent Key Words: division algebra, Schur index, ramification, Picard group, Brauer group\smallskip}

\footnote{}{}
\vfill\eject
\leftline{\bf Introduction}

This paper was first motivated by the following question 0.1, posed some
time ago by M.~Schacher, cf.~[Sc], and
resurrected more recently by L.~Small:

\proclaim Question 0.1. Suppose that $D_1$ and $D_2$ are division
algebras (finite over their centers), and
 both centers contain an algebraically closed field $F$. Does
$D_1 \otimes_F D_2$ have zero divisors?

We have an example where 0.1 is false, which appears in Section~4.

Let $F_i$ be the center of $D_i$ and write this as $D_i/F_i$. Then
it is well known (and we reprove below) that $F_1 \otimes_F F_2$
is a domain and we can set $K = q(F_1 \otimes_F F_2)$ to be its
field of fractions. Now Question~0.1 is equivalent to asking
whether $(D_1 \otimes_{F_1} K) \otimes_K (D_2 \otimes_{F_2} K)$ is
a division algebra. We will frequently switch between these two
points of view. It is enough to consider Question~0.1 when the
$D_i/F_i$ have prime power degree,
with respect to the same prime number. We will always assume this,
and frequently restrict our focus to the case where the $D_i$ have
(the same) prime degree. This assumption simplifies our arguments,
and still is challenging and interesting.

As work on this problem proceeded, it became clear that this should
be viewed as a piece of the following more general subject.
To understand Question~0.1 one needs to understand the center of
$D_1 \otimes_F D_2$ which leads to:

\proclaim Question 0.2. Suppose $F$ is algebraically closed and $F_i/F$
are field extensions. What are the properties of $F_1 \otimes_F F_2$?

Obviously Question~0.2 is ridiculously vague, but in this paper
we will ask and partially answer questions about the Picard groups
and Brauer groups of $F_1 \otimes_F F_2$. This seems most relevant
to Question~0.1. Moreover, Question~0.1 can be viewed
as being subsumed by Question~0.2 if we include in Question~0.2 the properties of
Azumaya algebras with center $F_1 \otimes_F F_2$.

In Question~0.1 and every version of Question~0.2 we consider here, we may reduce to the case
that the $D_i/F_i$ are finitely generated {\bf as division algebras}.
That is, the $F_i/F$ are finitely generated as
fields. In other language, we write $F_i = F(V_i)$ for a
projective variety $V_i$.

In Sections~2 through~4 below, we will assume that the ground
field $F$ has characteristic 0. This is to allow us to quote
resolution of singularities and write $F_i = F(V_i)$ where $V_i$
is a non-singular projective variety. In Section~3 we quote and
use the theorem of resolution of divisors.

Let us outline the paper to follow. In the rest of this introduction
we define some notation and observe one well known general
Brauer group fact. In Section~1 we make general observations
about Question~0.1, including the affirmative answer when
$D_1$ is commutative. Also in Section~1 is the perhaps
surprising connection between Question~0.1 and the ramification
of the $D_i$. More precisely, if $D_1$ is totally ramified at
a discrete valuation domain $R$ with $q(R) = F_1$, then
$D_1 \otimes_F D_2$ is a domain. In particular, when the $D_i$
have prime degree, and $D_1$ is ramified, then $D_1 \otimes_F D_2$
is a domain. This is partial justification for the idea that
``usually'' $D_1 \otimes_F D_2$ is a domain.

In Section~2 we cover some results about ramification that we need
in Section~3. The main idea is to show that we can eliminate all
ramification of a Brauer group element with a finite field
extension. In our case we need this extension to be a tensor
product of extensions of each of the $F_i = F(V_i)$. Though we do
not need it here, it is natural to ask that this extension be of
degree bounded in terms of the order of the Brauer group element
and the dimensions of the $V_i$. This stronger result is due to A.
Pirutka~[P] and we repeat a slight modification of her argument
here so we can further observe that the extension can be chosen a
product a fields as we need.

Section~3 contains the main body of our results about $F_1
\otimes_F F_2$.  Our feeling is that a full understanding of
Question~0.1 
in general requires a fuller understanding of these rings.
Finally, Section~4 has our example where $D_1 \otimes_F D_2$ is
not a domain, accompanied by the theory of Brauer groups over
curves that we need for the example. This material on Brauer
groups of curves has obvious independent interest.

The following fact is well known but a reference is hard to find.
If $V$ is a scheme, then 
$\Br(W)'$ for any $W \to V$ is defined to be the subgroup of the
Brauer group 
comprised of all elements of order prime to the characteristic of
all field points of $V$.

\proclaim Proposition 0.3. Let $V$ be an irreducible non-singular
scheme. Then $\Br(V)' = \cap_{P \subset V} \Br({\cal O}_{V,P})'$
where $P$ ranges over all irreducible codimension 1 subschemes and
${\cal O}_{V,P}$ is the stalk of $V$ at $P$, so ${\cal O}_{V,P}$
is a discrete valuation domain.

\proof When $V = \Spec(R)$ for a regular local ring $R$ this is
[Ho]. In [M, p.~147] it is observed that $U \to \Br(U)'$ is a Zariski
sheaf and the result follows.~\qed

\medskip
\leftline{\bf Section 1. First Observations}
\medskip

In this section we make some initial observations about question 0.1.
It is well known that if both $D_i = F_i$ are commutative then
$F_1 \otimes_F F_2$
is an integral domain. The next result is a generalization.

\proclaim Lemma 1.1. Suppose $D_1 = F_1$ is commutative. Then
$F_1 \otimes_F D_2$ has no zero divisors.

\proof Write $F_1 = F(V_1)$. Suppose $0 \not= \alpha_i = \sum_j
a_{i,j} \otimes b_{i,j} \in F_1 \otimes_F D_2$ are such that
$\alpha_1\alpha_2 = 0$. We can assume that for each $i$ the set of
$b_{i,j}$ are linearly independent over $F$. Since $V_1$ is
irreducible, there is an $F$ point on $V_1$ such that all the
$a_{i,j}$ are defined at this point and for each $i$ one of the
$a_{i,j}$ is nonzero. That is, there is a local ring $R \subset
F_1$ and $\phi: R \to F$ such that all $a_{i,j} \in R$, and for
each~$i$ some $\phi(a_{i,j}) \not= 0$. Then $\phi$ induces $\Phi:
R \otimes_F D_2 \to F \otimes_F D_2 = D_2$ such that
$\Phi(\alpha_i)$ is defined and $\Phi(\alpha_1)\Phi(\alpha_2) =
0$. Since the $b_{i,j}$ are linearly independent we have
$\Phi(\alpha_i) \not= 0$ for both $i$, and this is a
contradiction.~\qed

Note that this argument uses very little about $D_2$ except that
it is a domain. In particular, we need not assume that it is
finite over its center. Secondly, if $F_1$ is arbitrary (i.e.~not
necessarily finitely generated over its center) we can replace it
by the subfield generated by the $a_{i,j}$ and the same result
holds. All of this initially suggested to us that $D_1 \otimes_F
D_2$ is always a domain, but we have a counterexample. However,
our intuition still is the (not precise) feeling that $D_1
\otimes_F D_2$ is {\bf usually} a domain.

Though 1.1 is not hard, when combined with valuations,
it yields a result that says that $D_1 \otimes_F D_2$ very often is a
domain. Let $R$ be a discrete valuation domain with fraction field
$q(R) = F_1$. Then $R$ defines on (most of) the Brauer group
$\Br(F_1)$ a ramification map $\ram_R: \Br(F_1)' \to H^1(\bar
R,\Q/\Z)$. We say that $D_1$ is ramified at $R$ if $\ram_R([D_1])
\not= 0$. We say that $D_1$ is totally ramified if the order of
$\ram_R([D_1])$ is equal to the degree of $D_1$. Note that if
$D_1$ has prime degree, then it is ramified at $R$ if and only if
it is totally ramified at $R$.

To be unramified at all possible $R$ is a very strong condition.
We define the
unramified Brauer group to be the intersection of the kernels of
all these ramification maps, with respect to all these $R$.
The unramified Brauer group is much much smaller than the full
Brauer group. Thus the following result suggests that  $D_1
\otimes_F D_2$ is ``generically'' a domain.

\proclaim Theorem 1.2. a) Suppose $D_1$ is totally ramified
at some discrete valuation domain $R$ with $q(R) = F_1$.
Then $D_1 \otimes_F D_1$ is a domain.
\smallskip
b) Suppose $D_1/F_1$ has prime degree $p$ not equal to the
characteristic of $F$. If $D_1$ is ramified with respect to some
$R$ as in part (a), then $D_1 \otimes_F D_2$ is a domain.

\proof Part b) is a consequence of part~a) by our remark above. To
prove part a), let $\hat R$ be the completion and $\hat F = q(\hat
R)$. Denote by $\bar F_1 =R/M = \hat R/\hat M$ the residue field
of $R$ and $\hat R$. Our assumption on $D_1$ implies that $D_1
\otimes_F \hat F$ has degree equal to exponent and thus is a
division algebra. It follows that $R$ extends to a noncommutative
discrete valuation ring $S \subset D_1$ which defines a valuation
on $D_1$. More precisely, $S$ contains a unique maximal ideal
$S\pi$ such that $S\pi = \pi{S}$ is a two sided ideal, $D_1^* =
\cup_{n \in Z} S^*\pi^n = \cup_n \pi^nS^*$, and $S,\pi{S}$ lies
over $R,M$. Since $D_1$ is totally ramified, it follows that $L=
S/\pi{S}$ is a (commutative) field.

Suppose $0 \not= \alpha_i \in D_1 \otimes_F D_2$ are such that
$\alpha_1\alpha_2 = 0$. Write $\alpha_i = \sum a_{i,j} \otimes
b_{i,j}$ as above, where, again, for each $i$ the $b_{i,j}$ are
linearly independent over $F$. Then we can write  all $a_{1,j} =
\pi^{m_j}u_{1,j}$ and all $a_{2,j} = u_{2,j}\pi^{n_j}$ where all
$u_{i,j} \in S^*$. By changing $\alpha_1$ into $\pi^n\alpha_1$ for
some $n$ we can assume all $m_j \geq 0$ and some $m_j = 0$.
Similarly,  working on the other side, we can assume all $n_j \geq
0$ and some $n_j = 0$. In our other language, if $\phi: S \to L$
is the canonical morphism, we have all $\phi(a_{i,j})$ are defined
and for each $i$ there is $j$ such that $\phi(a_{i,j}) \not= 0$.
Again $\phi: S \to L$ induces $\Phi: S \otimes_F D_2 \to L
\otimes_F D_2$ and $\Phi(\alpha_i) = \sum_j \phi(a_{i,j})
\otimes_F b_{i,j} \in L \otimes_F D_2$. Since the $b_{i,j}$ are
linearly independent over $F$ the $1 \otimes b_{i,j} \in L
\otimes_F D_2$ must be linearly independent over $L$. In
particular, $\Phi(\alpha_i) \not= 0$ for both $i$. Since
$\Phi(\alpha_1)\Phi(\alpha_2) = 0$ we have a contradiction to
Lemma 1.1.~\qed

We remark that again 1.2 uses nothing about $D_2$ except that it
is a domain. We also remark that the totally ramified condition
can be eased a bit. Suppose $D_1$ is ramified at such an $R$ and $D_1 \otimes_F \hat F$ is a division algebra. Let $S$ exists as above
but $S/{\pi}S = \bar D_1$ is a division algebra with center
$L \not\supset \bar F_1$. Note that $\bar D_1/L$ has degree smaller
than $D_1$
and the proof of 1.2 shows that if $D_1 \otimes_F D_2$ has a zero divisor
then so does $\bar D_1 \otimes_F D_2$.

\medskip
\leftline{\bf Section 2. Ramification}
\medskip

In this section we investigate some questions about ramification
of Brauer group elements that we need in the rest of the paper. To
be precise, we need a result that any Brauer group element $\alpha
\in \Br(F(V_1) \otimes_F F(V_2))$ restricts to an everywhere
unramified element after an extension of the form $F(V_1' \times_F
V_2')$ where each $F(V_i')/F(V_i)$ is a finite field extension.
This can be done, but our original proof of this had the
uncomfortable property that the degree of this extension is
unbounded as we vary $\alpha$ among all elements of the same
order. This unboundedness did not constrain our arguments here,
but was unsatisfactory as it seemed that there should be a bound
on the degrees of the $F(V_i')/F(V_i)$ that only depends on the
dimension of the $V_i$ and the order of $\alpha$.

In fact, ignoring the 
specific requirements of this paper, the more natural question is
the following. Suppose $\alpha \in \Br(F(V))$. Is there a field
extension $F(V')/F(V)$ splitting all the ramification of~$\alpha$,
with degree bounded by a function of the order of $\alpha$ and the
dimension of $V$? 
It is believed that the {\bf index} of~$\alpha$ should have a
similar bound.  The result about splitting ramification would then
be evidence for this index conjecture. 
The second author asked the above splitting ramification question
at the workshop ``Deformation Theory, Patching, Quadratic Forms,
and the Brauer Group'' in January 2011 at the American Institute
of Mathematics. In April 2011 an affirmative answer was provided
by Alena Pirutka~[P]. Pirutka's result also generalizes to higher
degree cohomology. In this section we provide a slightly modified
proof of her result because we need to  observe further that we
can choose our $F(V')/F(V)$ to have  the form $F(V_1' \times_F
V_2')/F(V_1 \times V_2)$, and there is no reason to give our
earlier unbounded result. It should be noted that the bound in [P]
(and below) is known not to be strict even in the dimension~2
case.

To accomplish these results we need to make an observation about
what it takes to split all the ramification over a regular local
ring. Let $R$ be a regular local ring containing $F$, and take
$\alpha \in \Br(q(R))$ of order $q$. Suppose the ramification
locus of $\alpha$ has non-singular components with normal crossings
at $R$.

\proclaim Lemma 2.1. Suppose $\alpha$ ramifies at $\pi \in M -
M^2$ where $M$ is the maximal ideal of $R$. Set $S = R/\pi$ and
let $\bar L/q(S)$ be the ramification defined by $\alpha$. Then
all ramification of $\bar L/q(S)$ is at primes which are the
images of primes in the ramification locus of $\alpha$.

\proof Suppose not, and that in fact $\bar L$ ramifies at a prime
$\bar \delta$ of $S$ not on the list. Consider the inverse image of
$\bar \delta$ which is a height two prime $Q \subset R$. Note that
$Q$ contains none of the primes, except $\pi$, where $\alpha$ ramifies.
Let $T$ be the localization of $R$ at $Q$, so $T$ is a two dimensional
regular local ring. Set $\bar T = T/\pi$. Then $\pi$ is the only prime of $T$
where $\alpha$ ramifies. By e.g. [S, p.~129] this implies that
$\bar L/q(\bar T)$ is unramified at $\bar T$, a contradiction.~\qed

\proclaim Theorem 2.2. If $\alpha$ and $R$ are as above, and
$\pi_1,\ldots,\pi_r$ are the primes where $\alpha$ ramifies, then
$\alpha =  [\, \prod_{j=r}^1
(u_j\pi_1^{a_{1,j}}\ldots\pi_{j-1}^{a_{j-1,j}},\pi_j)_q]\alpha',$
where the $u_j$ are units and $\alpha' \in \Br(R)$.

\proof We induct on $r$. Let $\bar R = R/\pi_r$, and $\bar
L/q(\bar R)$ be the ramification of $\alpha$ at~$\pi_r$. By the
lemma $\bar L$ only ramifies on the images, $\bar \pi_i$, of the
$\pi_i$ for $i < r$. Since $\bar R$ is a~UFD, this implies $\bar L
= q(\bar R)((\bar u_r\bar
\pi_1^{a_{1,r}}\ldots\pi_{r-1}^{a_{r-1,r}})^{1/q})$ for a unit
$u_r$ and integers~$a_{i,r}$. Thus,
 $\alpha /
(u_r\pi_1^{a_{1,r}}\ldots\pi_{r-1}^{a_{r-1,r}},\pi_r)_q$ does not
ramify at $\pi_r$ and only ramifies at the $\pi_i$ for $i < r$. We
are done by induction on $r$.~\qed

We are going to kill ramification by the following trick.

\proclaim Proposition 2.3. Suppose $R$ is a regular local ring and
$\alpha$ and the $\pi_j$ are as above. Let $L \supset q(R)$ be a
field extension where for each $i$ there are units $v_i$ such that
$v_i^{-1}\pi_i$ is an $n$ power in $L$. Then $\alpha_L$ is
unramified with respect to any discrete valuation lying over a
localization of $R$.

\proof Write $\alpha$ as above and consider $\alpha'' =
\alpha'\prod_{j=1}^r (u_jv_1^{a_{1,j}}\ldots
v_{j-1}^{a_{j-1,j}},v_j)_q$. Then $\alpha$ and $\alpha''$ have the
same image in $\Br(L)$ and $\alpha'' \in \Br(R)$. Thus $\alpha''
\in \Br(R_P)$ for any prime $P$ and the result is clear.~\qed

Of course the difficulty is in constructing such an $L$ that works at
all the stalks. As above, let $V/F$ be smooth projective of dimension
$d$ and let $\alpha \in \Br(F(V))$. After blowing up (e.g. [K, p.138])
we may assume that the ramification locus of $\alpha$ consists of
non-singular irreducible components with normal crossings. Our next
result is really about such a set of divisors.

\proclaim Theorem 2.4. Let ${\cal D}$ be a set of non-singular
irreducible divisors of $V$ with normal crossings. Then there is a
morphism $V' \to V$ formed by blowing up along a succession of
non-singular subvarieties with the following property. Let ${\cal
D}'$ be the set of divisors of $V'$ formed as the union of strict
transforms of all elements of ${\cal D}$ and all exceptional
divisors (and the strict transforms of exceptional divisors). Then
${\cal D}'$ is the disjoint union of ${\cal D}_i'$ for $1 \leq i
\leq d$ such that each ${\cal D}_i'$ consists of disjoint
irreducible divisors.

\proof We will make repeated use of the following fact which uses
that the components of ${\cal D}$ are all non-singular with normal
crossings. Let $E$ be a component of a non-trivial intersection of
$D_1,\ldots, D_r$, all of which are elements of ${\cal D}$. If $V'
\to V$ is the blowup at $E$ and we identify $D_i$ with its strict
transform in $V'$, then in $V'$ the intersection of the $D_i$ is
empty.

Furthermore, any nonempty intersection of $r$ elements of ${\cal
D}$ has dimension $d-r$ and is the disjoint union of non-singular
components of that dimension. In particular, $r \leq d$. First we
look at all the nonempty intersections of $d$ elements of ${\cal
D}$, which altogether comprise a finite set of points. We form
$V_1 \to V$ by blowing up at all those points. Let ${\cal D}_1' =
{\cal D} \cup {\cal D}_1$, where ${\cal D}$ are the strict
transforms in $V_1$ of the divisors ${\cal D}$ in $V$, and ${\cal
D}_1$ are all the exceptional divisors which are obviously all
disjoint.

In ${\cal D}$ (viewed as divisors in $V_1$) we 
define  $E_1$ to be all nonempty intersections of subsets of order
$d-1$. Since any $d$ elements of ${\cal D}$ have
empty intersection,  the components of $E_1$ are all disjoint,
non-singular curves. We let $V_2 \to V_1$ be the blowup at all
these curves and set ${\cal D}_2' = {\cal D}_2 \cup {\cal D}_1
\cup {\cal D}$ where ${\cal D}_2$ are the exceptional divisors and
the rest of the terms again are strict transforms. Proceeding in
this way we are done because at the last step all the elements of
${\cal D}$ will be disjoint.~\qed

The above argument says that at the level of divisors we can
separate the ramification locus of $\alpha \in \Br(F(V))$ so that
if all the $\sum_{E \in {\cal D}_i} E$ were principal, we could
take all these $q$ roots and kill all ramification. Since these
divisors need not be principal, we have to proceed as follows.

In the arguments to come, we will be given a finite set of
irreducible closed sets ${\cal C}$ of a variety $V'$.
Let ${\cal C'}$ be the union of ${\cal C}$ and the finite set of all components
of all intersections of subsets of ${\cal C}$.
Of course, ${\cal C'}$ is a finite set closed under the process of taking
components of intersections of subsets.
Let ${\cal M}$
be the set of minimal elements of ${\cal C'}$, being all elements
which do not properly contain another element of ${\cal C'}$. Then
all the elements of ${\cal M}$ are disjoint. Thus we can take the
stalk of the structure sheaf ${\cal O}_V$ at ${\cal M}$ and by abuse of
terminology we call this the stalk of ${\cal C}$.

\proclaim Lemma 2.5. Let $R$ be the stalk of ${\cal C}$. Then $R$
is a semilocal domain. If $V'$ is non-singular, then $R$ is
regular and a UFD. The prime elements (up to units) of~$R$
correspond to all irreducible divisors of $V$ which contain a
component of an intersection of a subset of ${\cal C}$.

\proof That $R$ is a UFD is well known and can be found,
for example, in [S1, p.~1546]. This rest is all clear.~\qed

We return to a set of irreducible divisors ${\cal D}$ as in the
conclusion of 2.4. That is, the elements of ${\cal D}$ are all
non-singular and together they have normal crossings. Moreover,
${\cal D} = \cup_{i=1}^d {\cal D}_i$ where the irreducible
divisors in each ${\cal D}_i$ are all disjoint. Let ${\cal D}_i$
consist of divisors $D_{i,j}$ where $1 \leq j \leq s(i)$. Let
$R_1$ be the stalk of $V$ at ${\cal D}$ which for the purposes of
this argument we rename ${\cal E}_1$. Then all the $D_{i,j}$
induce primes $\pi_{i,j,1}$ on $R_1$ and we can choose $f_{i,1}$
such that $f_{i,1} R_1 = (\prod_j \pi_{i,j,1})R_1$. Looking
globally, the principal divisor $(f_{i,1})$ equals $\sum_j D_{i,j} +
\sum_k n_{i,k,1}E_{i,k,1}$, where no component of any of the intersections
of any subset of the elements of ${\cal E}_1$ is contained in any
$E_{i,k,1}$.

By induction, assume $R_l$, ${\cal E}_l$,
$f_{i,l}$, $E_{i,k,l}$ and $n_{i,k,l}$ have been defined for all
$l < m$ where:

a) ${\cal E}_l$ is the set of all $D_{i,j}$ and all $E_{i,k,l'}$
for all $l' < l$.

b)  $(f_{i,l}) = \sum_j D_{i,j} + \sum_k n_{i,k,l}E_{i,k,l}$.

c) $R_l$ is the stalk of ${\cal E}_l$.

d) No $E_{i,k,l}$ contains a component of an intersection
of elements of ${\cal E}_l$.

Of course we define
${\cal E}_m$
to be the set of all $D_{i,j}$ and all $E_{i,k,l}$
for all $i,k$ and $l < m$. Equally obviously, we set $R_m$ to be
the stalk at
${\cal E}_m$ and in $R_m$ we let $D_{i,j}$ define $\pi_{i,j,m}$ on
$R_m$. Let $f_{i,m}$ be such that $f_{i,m}R_m = (\prod_j
\pi_{i,j,m})R_m$. Of course, we define the $E_{i,k,m}$ and $n_{i,k,m}$ via 
$(f_{i,m}) = \sum_j D_{i,j} + \sum_k n_{i,k,m}E_{i,k,m}$.

We perform the above construction until $m = d$ where $d$ is the
dimension of~$V$. We claim that:

\proclaim Lemma 2.6.
Let $C \subset V$ be an irreducible closed subset contained
in some $D_{i,j}$.
Then for some $m$, $C$ is not contained in $E_{i,k,m}$ for any
$i$ and $k$.

\proof Otherwise, for each $m$
there are $i(m), k(m)$ such that $C \subset E_{i(m),k(m),m}$. Now
for each $m$ no component of $D_{i,j} \cap E_{i(1),k(1),1} \cap
\ldots \cap E_{i(m-1),k(m-1),m-1}$ is contained in
$E_{i(m),k(m),m}$. It follows that for every $m$ every component
of $D_{i,j} \cap E_{i(1),k(1),1} \cap \ldots \cap E_{i(m),k(m),m}$
has dimension less than or equal to $d - m - 1$. When $m = d$ this
is a contradiction.~\qed

Now assume $\alpha \in \Br(F(V))$ has exponent $q$ and ${\cal D}$
is the set of divisors where $\alpha$ ramifies. Assume we have
blown up so that all the elements of ${\cal D}$ are non-singular
with normal crossings and further that ${\cal D}$ is the union
of ${\cal D}_i$ as in Theorem~2.4. Note that it is possible during the process
of blowing up that an exceptional divisor will not be a divisor
where $\alpha$ ramifies. If we exclude it from ${\cal D}$ and
${\cal D}_i$ the conclusions of Theorem~2.4 still stand. Next form
the $f_{i,m} \in F(V)$ as in Lemma~2.6 for $1 \leq m \leq d$.

Let $F'/F$ be defined by taking the $q$ roots of all $f_{i,m}$.
Note that $F'/F$ has degree less than or equal to $q^{d^2}$. We
see next that $F'$ splits all ramification of~$\alpha$. Note that
in [S1, p.~1584] we proved that when $q$ is prime, and $S$ is a
non-singular surface, an extension of degree $q^2$ and not $q^4$
splits all ramification. It is therefore conceivable that a more
careful study of ramification in dimensions greater than~2 would
yield a better bound than the one below.

\proclaim Theorem 2.7 (see [P]). The restriction $\alpha|_{F'} \in \Br(F')$
is unramified everywhere and the degree of $F'/F$ is bounded by a function
of $d$ and $q$.

\proof The second statement is clear. As for the first, suppose
$S$ is a discrete valuation domain with $q(S) = F'$. Then $S$ lies
over an irreducible
closed subset $C \subset V$. Let $R$ be the stalk ${\cal O}_{V,C}$
of $V$ at $C$. If $C$ is not contained in any $D_{i,j}$ then
$\alpha \in \Br(R_C)$ and so the restriction of $\alpha$ is
unramified at $S$.

Thus we assume $C \subset D_{i,j}$ for some $i,j$.  Note that by
disjointness for each~$i$ there is at most one such $j$. Let $I$
be the set of $i$ where $C \subset D_{i,j}$ for some $j$.  By~2.6
we can take $m$ such that $C$ is not contained in any $E_{i,k,m}$.
The $f_{i,m}$, for $i \in I$, must be prime elements of $R_C$. We
are done by Proposition~2.3.~\qed

\proclaim Theorem 2.8. Suppose $\alpha \in \Br(F_1 \otimes_F F_2)$
has exponent $q$.
Then there are finite field extensions $F_i' \supset F_i$ such that
$\alpha$ maps to an everywhere unramified element of
$\Br(F_1' \otimes_F F_2')$.

\proof Since $\alpha \in \Br(F_1 \otimes_F F_2)$, it follows that
the ramification locus of~$\alpha$ on $V_1 \times_F V_2$ consists
of vertical and horizontal irreducible divisors, where
``vertical'' means the divisor has the form $\pi_1^*(D)$ for $D
\subset V_1$ a divisor ( and ``horizontal'' is the $V_2$ version).
We can blow up $V_1$ so that the vertical irreducible divisors are
non-singular and have normal crossings. We do the same thing to
$V_2$. Of course, this implies that their respective pullbacks,
taken all together, have non-singular components and normal
crossings.

Let ${\cal D}_i$ be the irreducible divisors in the ramification
locus of $\alpha$ coming from~$V_i$. Let $d_i$ be the dimension of
$V_i$. By further blowing up we may assume that  each ${\cal D}_i$
is the union of disjoint subsets ${\cal D}_{i,j}$ such that the
elements of ${\cal D}_{i,j}$ are themselves disjoint. Viewing the
${\cal D}_i$ as divisors on $V_i$, form $f_{i,j,m} \in F(V_i)$ as
in Lemma 2.6. Let $F_i' \supset F(V_i)$ be the field obtained by
adjoining $f_{i,j,m}^{1/q}$ for all $j$ and $1
\leq m \leq d_i$. We can write $F_i' = F(V_i')$ and $F' =
F(V_1' \times_F V_2')$.

Let $\alpha'$ be the restriction of $\alpha$ to $\Br(F')$. By 0.3
it suffices to show that $\alpha'$ is unramified with respect to
the stalk of any irreducible divisor on $V_1' \times_F V_2'$. If
such a divisor is not horizontal or vertical, then $S$ contains
$F_1' \otimes_F F_2'$ and $\alpha' \in \Br(F_1' \otimes_F F_2')$
since $\alpha \in \Br(F_1 \otimes_F F_2)$. Thus we may assume by
symmetry that $S$ is the stalk at a vertical divisor $D \times
V_2'$. That is, $S$ lies over an irreducible $C \times V_2 \subset
V_1 \times V_2$. If $R$ is the stalk ${\cal O}_{V_1 \times V_2,
C\times V_2}$ then, in the ramification locus of~$\alpha$, only
vertical primes appear as primes in $R$. By 2.2, $\alpha$ is a
product of symbols involving vertical primes and an element of
$\Br(R_C)$. Thus by the argument of 2.3 and 2.7, if we restrict
$\alpha$ to $\Br(L)$ where $L = q(F_1' \otimes_F F_2)$, then
$\alpha|_L$ is unramified at any discrete valuation lying over $C
\times V_2$.~\qed

\medskip
\leftline{\bf Section 3. Tensor Products of Fields}
\medskip

In the previous section we saw that the tensor product of fields
(over an algebraically closed field) is always a domain. In that sense
this is not a case we need to consider. But it will be useful to us,
and of considerable interest, to further study tensor products of fields.
After all, these rings are the centers of the tensor products of division
algebras, and therefore the arithmetic of these rings is important to the study
of the more general tensor products.

To begin, in this section $F$ is always an algebraically closed field of
characteristic 0. We make this characteristic assumption because
we make frequent use of the fact that varieties over $F$ have resolutions
of singularities and resolutions of divisors.

Let $F_1$ and $F_2$ be fields of finite
type over $F$. That is, the $F_i$ are the fields of fractions
of projective $F$ varieties. Because of our assumptions each $F_i$ is,
in fact, the function field of a
smooth projective variety $V_i$ defined over $F$. Set $R = F_1
\otimes_F F_2$. Our goal in this section is to study the properties
of $R$ and related rings. In particular, we will be interested in
the Picard group of $R$ and the Brauer group of $R$. Let $\bar
F_i$ denote the algebraic closure of $F_i$. Frequently we will be
extending the scalars of $V_i$ to general $K \supset F$ and more specifically
to $\bar F_j \supset F_j \supset F$. We
write $V_i \times_F K$ as $V_i/K$ and similarly for $\bar F_j$. We
write $\Pic(V_i/K)$ for the Picard group of $V_i/K$ and ${\cal
P}{\it{ic}}(V_i/K)$ for the Picard scheme defined over $K$. As a
source for the basic properties of this scheme one can
use [BLR, p.~199--235].

By [BLR, p.~232 and p.~210] the connected component ${\cal
P}{\it{ic}}^0(V_2/K)$ is a projective scheme over $K$ which is of
finite type. By [BLR, p.~231] it is smooth and we call it the
Picard variety of $V_2/K$. Being a group scheme, ${\cal
P}{\it{ic}}^0(V_2/F)$ is an abelian variety. Moreover,
$\Pic(V_2/K)$ can be identified with the $K$ points of ${\cal
P}{\it{ic}}(V_2)$ ([BLR, p.~204]). It therefore makes sense to let
$\Pic^0(V_2/K)$ be the $K$ points of ${\cal P}{\it{ic}}^0(V_2/F)$.
Also, it follows that $\Pic(V_2/F) \to \Pic(V_2/K)$ is injective
for any field $K \supset F$. Note that ${\cal P}{\it{ic}}(V_2/F)
\times_F K = {\cal P}{\it{ic}}(V_2/K)$ because of the functorial
definition of ${\cal P}{\it{ic}}(V_2/K)$. Because irreducibles
over $F$ are absolutely irreducible it follows that ${\cal
P}{\it{ic}}^0(V_2/F) \times_F K = {\cal P}{\it{ic}}^0(V_2/K)$.

\proclaim Proposition 3.1. Let $K$ be a field containing $F$. Then
$$\Pic(K \otimes_F F_2) = \Pic(V_2/K)/\Pic(V_2/F) =
\Pic^0(V_2/K)/\Pic^0(V_2/F).$$

\proof $\Pic(K \otimes_F F_2)$ is the direct limit of the $\Pic(U/K)$
where $U \subset V_2$ are open subvarieties defined over $F$. Thus
$\Pic(V_2/K) \to \Pic(K \otimes_F F_2)$ is clearly surjective. If
some $\alpha \in \Pic(V_2/K)$ maps to 0 in $\Pic(U/K)$, then
lifting to divisors we have for some $f \in K(V_2)^*$ that $\alpha
= (f) + \sum n_iD_i$ where $V_2 - U$ is the union of irreducibles
$D_i$ defined over $F$. That is, $\alpha$ is in the image of
$\Pic(V_2/F)$. This proves the first equality.

Since the irreducible components of ${\cal P}{\it{ic}}(V_2/F)$
stay irreducible over $K$, it follows that $\Pic^0(V_2/K) \to
\Pic(V_2/K)/\Pic(V_2/F)$ is surjective and the kernel is
$\Pic(V_2/F) \cap \Pic^0(V_2/K) = \Pic^0(V_2/F)$.~\qed

As the torsion subgroup of $\Pic^0(V_2/K)$ is all defined over $F$
and $\Pic^0(V_2/F)$ is divisible we have:

\proclaim Corollary 3.2. $\Pic(K \otimes_F F_2)$ is torsion free.

\proof If $\alpha \in \Pic^0(V_2/K)$ satisfies $n\alpha \in \Pic^0(V_2/F)$
then $n\alpha = n\beta$ for some $\beta \in \Pic(V_2/F)$.
Thus $n(\alpha - \beta) = 0$ implying $\alpha - \beta \in \Pic(V_2/F)$
and so $\alpha \in \Pic(V_2/F)$.~\qed

Another immediate corollary of Proposition 3.1.
is:

\proclaim Corollary 3.3.   $\bar F_1 \otimes_F F_2$ has divisible
Picard group.

\proof This follows because $\Pic(\bar F_1 \times_F F_2) =
\Pic^0(V_2/\bar F_1)/\Pic^0(V_2/F)$ and\break $\Pic^0(V_2/\bar F_1)$ is a divisible
group.~\qed

To understand the Brauer group of $F_1 \otimes_F F_2$, we begin by
showing that $\bar F_1 \otimes_F \bar F_2$ has Brauer group 0. The
first step is the unramified case.

\proclaim Lemma 3.4. Any element in the Brauer group of $V_1 \times V_2$
maps to 0 in $\Br(\bar F_1 \otimes \bar F_2)$.

\proof If $\alpha \in \Br(V_1 \times_F V_2)$, then certainly
$\alpha \in \Br(V_2/F_1)$. For any $n > 0$, let $\mu_n \subset
\mu$ be the subgroup of roots of 1 of order $n$. We have (e.g. [M,
p.~224]) $H^2(V_2/\bar F_1,\mu_n) = H^2(V_2/F,\mu_n)$. The Kummer
sequence induces the commutative diagram:
$$\matrix{
0&\to&\Pic(V_2/F)/n\Pic(V_2/F)&\longrightarrow&H^2(V_2/F,\mu_n)&\longrightarrow&_n\Br(V_2/F)&\to&0\cr
&&\downarrow&&||&&\downarrow\cr
0&\to&\Pic(V_2/\bar F_1)/n\Pic(V_2/\bar F_1)&\longrightarrow&H^2(V_2/\bar F_1,\mu_n)&\longrightarrow&_n\Br(V_2/\bar F_1)&\to&0\cr}$$
where $_n\Br$ refers to the $n$ torsion.
Also we know that $\Pic(V_2/F)/n\Pic(V_2/F) =
\Pic(V_2/\bar F_1)/n\Pic(V_2/\bar F_1)$.
Then applying the above diagram for all $n$ we have that
$\Br(V_2/F) \to \Br(V_2/\bar F_1)$ is an isomorphism.
In particular, $\alpha = \alpha_1 +
\alpha_2$ where $\alpha_1$ is in the image of $\Br(V_2/F)$ and
$\alpha_2$ maps to~0 in $\Br(V_2/\bar F_1)$. Certainly $\alpha_1$
maps to 0 in $\Br(\bar F_2)$ and so both $\alpha_i$ map to 0
in~$\Br(\bar F_1 \otimes_F \bar F_2)$.~\qed

By combining Lemma~3.4 with Theorem~2.8 we get:

\proclaim Theorem 3.5. $\Br(\bar F_1 \otimes_F \bar F_2) = 0$.

\proof Any $\bar \alpha \in \Br(\bar F_1 \otimes_F \bar F_2)$ is
the image of some $\alpha \in \Br(F(V_1) \otimes_F F(V_2))$. By
Theorem~2.8, we may assume that $\alpha$ is in the Brauer group of
$V_1 \times_F V_2$, and so we are done by Lemma~3.4.~\qed

As a consequence of the above theorem, any element of
$\Br(F_1 \otimes_F F_2)$ is split by an extension $F_1' \otimes_F F_2'$
where the $F_i'/F_i$ are Galois with group $G_i$. That is, $F_1'
\otimes_F F_2'$ is Galois over $F_1 \otimes_F F_2$ with group $G_1
\oplus G_2$.

Let $S = F_1' \otimes_F F_2'$ and $R = F_1 \otimes_F F_2$. From
[DI, p.116] we know that there is an exact sequence $\Pic(R) \to
\Pic(S)^G \to H^2(G,S^*) \to \Br(S/R) \to H^1(G,\Pic(S))$ where $G
= G_1 \oplus G_2$.

We next will show that by extending $S$ we may assume that  Brauer
group elements are crossed products. Let $\bar G_i$ be the Galois group of
$\bar F_i/F_i$ and $\bar S = \bar F_1 \otimes_F \bar F_2$.

\proclaim Lemma 3.6. a) Suppose that in the above sequence $\alpha \in \Br(S/R)$.
Then there are fields $F_i'' \supset F_i'$ such that $F_i''/F_i$ is Galois with
group $G_i'$ mapping to $G_i$ and under restriction $\alpha$ is in the
image of $H^2(G_1' \ \oplus G_2', S'^*)$ where $S' = F_1'' \otimes_F F_2''$.
\smallskip
b) There is a surjection $H^2(\bar G_1 \oplus \bar G_2,\bar S^*) \to \Br(R)$.

\proof Part b) is a consequence of a).
If $\beta \in H^1(G_1 \oplus G_2,\Pic(S))$ it suffices to
show that there are such $G_1', G_2',$ and  $S' = F_1'' \otimes_F F_2''$ with
$\beta$ mapping to $0$ in $H^1(G_1' \oplus G_2',\Pic(S'))$.
Now $H^1(\bar G_1 \oplus \bar G_2, \Pic(\bar S))$ is the direct limit
of all such $H^1(G_1' \oplus G_2', \Pic(S'))$ and it suffices to show
that $H^1(\bar G_1 \oplus \bar G_2, \Pic(\bar S)) = 0$.

By similar reasoning $H^1(\bar G_1 \oplus \bar G_2, \Pic(\bar S))$
is the direct limit of all $$H^1(\bar G_1 \oplus G_2', \Pic(\bar
F_1 \otimes_F F_2''))$$ taken over all $G_2'$ Galois extensions
$F_2''/F_2$. But by Corollaries 3.2 and 3.3, $\Pic(\bar F_1
\otimes_F F_2'')$ is torsion free divisible, so $H^1(\bar G_1
\oplus G_2', \Pic(\bar F_1 \otimes_F F_2'')) = 0$.~\qed

Thus to describe $\Br(F_1 \otimes_F F_2)$, we need to describe
$H^2(\bar G_1 \oplus \bar G_2, \bar S^*)$. To this end we next
observe:

\proclaim Lemma 3.7. a) $S^* \cong (F_1^* \oplus F_2^*)/F^*$ and
$\bar S^* = (\bar F_1^* \oplus \bar F_2^*)/F^*$.
\smallskip
b) $H^2(\bar G_1 \oplus \bar G_2, \bar S^*) \cong H^2(\bar G_1
\oplus \bar G_2, F^*) \cong H^2(\bar G_1 \oplus \bar G_2, \mu)$
where $\mu \subset F^*$ is the group of roots of 1, and so $\mu
\cong \Q/\Z$ as a $\bar G_1 \oplus \bar G_2$ module.

\proof We begin with a). The second statement of a) follows from
the first. Suppose $u \in (F_1 \otimes_F F_2)^*$ and let $F_i =
F(V_i)$ with $V_i$ projective non-singular. If we consider the
principle divisor $(u)$ of $u$ on $V_1 \times V_2$, then all the
zeroes and pole components must be horizontal or vertical. Let $D$
be the divisor of vertical zeroes and poles, which we can also
view as a divisor of $V_1/F$. Thus in $V_1/F_2$, $D$ is a
principal divisor, and since $\Pic(V_1/F) \to \Pic(V_1/F_2)$ is
injective, we know that $D$ is principal as a divisor over
$V_1/F$. In other words, there is an element $v \in F(V_1)^* =
F_1^*$ such that $u/v$ is a unit on $V_1/F_2$. In other words, $u
= vw$ where $v \in F_1^*$ and $w \in F_2^*$. On the other hand, if
$vw = 1$, then $v,w \in F_1^* \cap F_2^* = F^*$.

Turning to b), there is an exact sequence $$0 \to F^* \to (\bar
F_1^* \oplus \bar F_2^*/F^*) \to (\bar F_1^*/F^* \oplus \bar
F_2^*/F^*) \to 0,$$ and each $\bar F_i^*/F^*$ is torsion free
divisible.  2) follows immediately.~\qed

Let $\Sym(\bar G_1,\bar G_2)$ be defined as the direct limit of
$\Hom(\bar G_1,\mu_n) \otimes_{\Z} \Hom(\bar G_2,\mu_n)$ over all $n$.
Now we can invoke standard group cohomology and observe:

\proclaim Lemma 3.8. $H^2(\bar G_1 \oplus \bar G_2, \mu) 
 \cong
H^2(\bar G_1,\mu) \oplus \Sym(\bar G_1,\bar G_2) \oplus H^2(\bar
G_2,\mu)$.

\proof This follows, for example, from the Hochschild-Serre
spectral sequence applied to $\bar G_1 \oplus \bar G_2 \to \bar
G_1$. From the product structure, $H^2(\bar G_i,\mu) \to H^2(\bar
G_1 \oplus \bar G_2)$ is injective and so $H^2(\bar G_2,\mu)^{\bar
G_1} = H^2(\bar G_2,\mu)$ survives unchanged in the limit of the
spectral sequence.

It suffices to show that
$$H^1(\bar G_1,H^1(\bar G_2,\mu)) =  \Hom(\bar G_1,\Hom(\bar G_2,\mu)) =
\Sym(\bar G_1,\bar G_2)),$$ and that all the elements of
$\Sym(\bar G_1,\bar G_2)$ survive in the limit. The first
statement follows because $\Hom(\bar G_1,\Hom(\bar G_2,\mu))$ is
the direct limit of the $$\Hom(\bar G_1,\Hom(\bar G_2,\mu_n)) =
\Hom(\bar G_1,\mu_n) \otimes_{\Z} \Hom(\bar G_2,\mu_n)$$ for all
$n$. The second statement follows because all the elements of
$\Sym(\bar G_1,\bar G_2)$ are images of cup products of elements
of the $H^1(\bar G_i,\mu_n)$.~\qed

If, as above, the $\bar G_i$ are absolute Galois groups of the
fields $F_i$, we write $$\Sym(\bar G_1,\bar G_2) =
\Sym(F_1,F_2).$$ We can think of this last group as abstract
symbols in the cohomology. We are ready for:

\proclaim Theorem 3.9. $\Br(F_1 \otimes_F F_2) = \Br(F_1) \oplus
\Br(F_2) \oplus I$, where $I$ is an image of $\Sym(F_1,F_2)$.

\proof By Lemma 3.6 and 3.7 there is a surjection
$$\phi: H^2(\bar G_1,\mu) \oplus \Sym(\bar G_1,\bar G_2) \oplus H^2(\bar G_2,\mu) =
H^2(\bar G_1 \oplus \bar G_2, \mu) \to \Br(F_1 \otimes_F F_2).$$
We can identify $H^2(\bar G_i,\mu)$ with $\Br(F_i)$. Since $F_i =
F(V_i)$ and $V_i$ has an $F$ point, the induced map $\Br(F_i) \to
\Br(F_1 \otimes_F F_2)$ is injective.~\qed

Note that $\Br(F_1)$ and
$\Sym(F_1,F_2)$ map to zero in $\Br(\bar F_1 \otimes_F F_2)$ and so
taking direct limits we have:

\proclaim Lemma 3.10. $\Br(\bar F_1 \otimes_F F_2) = \Br(F_2)$
and $\Br(F_1 \otimes_F \bar F_2) = \Br(F_1)$.

Now we can finish 3.9 by noting 3.10 and the restrictions
$$\Br(F_1 \otimes_F F_2) \to \Br(\bar F_1 \otimes_F F_2)\qquad \Br(F_1 \otimes_F F_2) \to \Br(F_1 \otimes_F \bar
F_2)$$
 imply that any element in the kernel of $\phi$
has to be in $\Sym(\bar G_1,\bar G_2)$.~\qed

Of course $\Sym(F_1,F_2) \to \Br(F_1 \otimes_F F_2)$ can be
interpreted as the union of symbol algebra induced maps $(a_1,a_2)
\to (a_1,a_2)_n$ where $a_i \in F_i^*$ but the symbol algebra has
center $F_1 \otimes_F F_2$. Also, we know from the example in
Section~4 this map is not injective.
In fact, noticing the
non-injectivity here was the idea that spurred the discovery of
the example of section 4. The connection may not be clear, but one way of viewing
questions about Schur index over $q(F_1 \otimes_F F_2)$ is that one
is considering the ``smallest'' way of representing a Brauer group
element. The non-injectivity above raised the possibility
of writing Brauer group elements in terms of fewer symbols
by ``using'' trivial elements in the image of $\Sym(F_1,F_2)$.
The connection is perhaps not rigorous, but it was strong enough
to suggest the example in the next section.

\medskip
\leftline{\bf Section 4. Products of Curves and the Counterexample}
\medskip

In this section we consider Brauer groups over products of curves
and use that machinery to provide a counterexample to our main
question. Let $F$ be a field of characteristic $0$. Although $F$
is not assumed to be algebraically closed,   it should be clear
that the algebraically closed case provides us important examples.

Suppose $C$ and $C'$ are two curves defined over $F$ with the
additional property that all torsion points of the Jacobians
$\Jac(C)$ and $\Jac(C')$ are $F$ rational and both curves have $F$
rational points. Let $K = F(C')$ and let $\bar K$ be the algebraic
closure of $K$. Let $\bar C = C \times_F \bar K$, and let $G$ be
the Galois group of $\bar K/K$. Let $\Tor(\Jac(\bar C))$ be the
torsion subgroup which by assumption has trivial $G$ action. Since
$\Jac(\bar C)/\Tor(\Jac(\bar C))$ is torsion free divisible, every
element of $H^1(G,\Jac(\bar C))$ is in the image of
$H^1(G,\Tor(\Jac(\bar C))) = \Hom(G,\Tor(\Jac(\bar C)))$. We have
the following three exact sequences of $G$ modules associated to
the curve $C$.

$$0 \to \bar K(C)^*/\bar K^* \to \Div(\bar C) \to \Pic(\bar C) \to 0$$
and
$$0 \to \bar K^* \to \bar K(C)^* \to \bar K(C)^*/\bar K^* \to 0$$
and
$$0 \to \Jac(\bar C) \to \Pic(\bar C) \to \Z \to 0$$
We will frequently apply the long exact cohomology sequence
to each of these sequences.

Since $C$ has a $K$ rational point, the last sequence splits.
Thus, $\Pic(\bar C)^G \to \Z$ is surjective, and since $H^1(G,\Z)
= 0$ we have that $H^1(G,\Jac(\bar C)) = H^1(G,\Pic(\bar C))$.
There is a discrete valuation ring $R \subset K(C)$ with $q(R) =
K(C)$ and residue field $\bar R = K$. Thus, $\bar R = \bar K
\otimes_K R$ is a discrete valuation domain with $q(\bar R) = \bar
K(C)$, and $\bar R \to \bar K$ is a $G$ map. Now $\bar K(C)^* =
\bar R^* \oplus \Z$, so there is a $G$ morphism $\bar K(C)^* \to
\bar K^*$ and the second sequence splits. Thus,
$$H^2(G,\bar K(C)^*/\bar K^*) = H^2(G,\bar K(C)^*)/H^2(G,\bar
K^*).$$ By Tsen's Theorem (e.g.~[Se, p.~162]) $$H^2(G,\bar K^*(C))
= \Br(\bar K(C)/K(C)) = \Br(K(C))$$ so $$H^2(G,\bar K(C)^*)/\bar
K^* = \Br(K(C))/\Br(K).$$ For any point $P$ of $\bar C$ let $G_P$
be the stabilizer in $G$ of $P$. Then $$\Div(\bar C) = \oplus_P \
\Z[G/G_P],$$ the direct sum being over all $G$-orbits of points.
Thus $$H^1(G,\Div(\bar C)) =\oplus_P \  H^1(G_P,\Z) = 0.$$
Consider the composite 
 $\phi: H^2(G,\bar K(C)^*) \to
H^2(G,\Div(\bar C)) = \oplus_P H^1(G_P,\Q/\Z)$.
 Since $G_P$ is the
absolute Galois group of the residue field of the $C$ point
defined by $P$, it is easy to see that $\phi$ is the sum of all
the ramification maps at all points of $C$. It follows that
$H^1(G,\Jac(\bar C)) = H^1(G,\Pic(\bar C)) \subset
\Br(K(C))/\Br(K)$ is the subgroup unramified at all of the points
of $C$, or in different language:

\proclaim Lemma 4.1. $H^1(G,\Jac(\bar C)) = \Br(C)/\Br(K)$.

Note that most of the paragraph preceding Lemma~4.1 is essentially due to Roquette ([Ro]).

Recall that we are interested in an element $\alpha \in
H^1(G,\Jac(\bar C))$ which is the image of an element $\alpha' \in
\Hom(G, \Tor(\Jac(\bar C)))$ with cyclic image of order $n$. If
$H$ is the kernel of $\alpha'$, let $L = \bar K^H$. If $\sigma{H}$
is a generator of $G/H$, let $P = \alpha'(\sigma)$ be the
associated element of order $n$ of $\Jac(\bar C)$ and let $P' \in
\Div(\bar C)$ be a preimage of~$P$. Let $\beta$ be the image of
$\alpha$ in $H^2(G,\bar K(C)^*/\bar K^*)$ under the coboundary.

We want to make $\beta$ more explicit. Since $\alpha$ is the image
of
some cocycle $\alpha_c \in H^1(G/H,\Jac(\bar C)^H)$, $\beta$ is
the image of $\beta_c \in H^2(G/H,(\bar K(C)^*/\bar K)^H)$ where
$\beta_c$ is the image of $\alpha_c$ under the $G/H$ coboundary. 
Since $G/H$ is cyclic, if $M$ is any $G/H$ module, $H^2(G/H,M) \cong 
M^G/N_{G/H}(M)$ where $N_{G/H}: M \to M$ is the norm map and the 
isomorphism depends on the choice of generator $\sigma{H}$ of $G/H$. 
In our case $M^G = (\bar K(C)^*/\bar K^*)^G = K(C)^*/K^*$. 
Tracing through the $G/H$ coboundary map we see that $\beta_c$
corresponds to the image of $fK^*$ where $f \in K(C)^*$ is such that the
divisor $(f) = nP'$. That is, as a Brauer group element $\alpha$
maps to the cyclic algebra $\Delta(L(C)/K(C),\sigma,f)$ (modulo
$\Br(K)$). All together we have:

\proclaim Lemma 4.2. The element $\beta$, viewed as an element of
$\Br(K(C))/\Br(K)$, is represented by the cyclic algebra
$\Delta(L(C)/K(C),\sigma,f)$ where $L/K$, $\sigma$ and $f$ are as
above.

We are interested in when $\beta$ is trivial. That is, given
$\alpha': G \to \Tor(\Jac(\bar C))$ as above, we are interested
when it maps to $0$. Let $\sigma$ be as above.

\proclaim Lemma 4.3. The element $\alpha'$ maps to 0 in
$H^1(G,\Jac(\bar C))$ if and only if there is an $L$ point $Q$ of
$\Jac(C)$ with $\sigma(Q) - Q = P$.

\proof From the definition of degree 1 cohomology, there is an $H$
fixed $Q' $ in~$\Jac(\bar C)$ such that $\sigma(Q') - Q' = P$. The
map $H^2(H,\bar K) \to H^2(H,\bar K(C))$ has been identified with
$\Br(L) \to \Br(L(C))$, which is injective since $C$ has a $K$
rational point. Thus $H^1(H,\bar K(C)^*/\bar K^*) = 0$ and it
follows from the long exact cohomology sequence that $Q'$ is the
image of an $H$ fixed element $Q''$ of the divisor group
$\Div(\bar C)$. That is, as an element of $\Pic(\bar C)$, $Q'$ can
be written as a sum of $H$ orbits of points. In other language,
$Q'$ corresponds to an $L$ point, $Q$, of $\Pic(C)$. After
subtracting a suitable multiple of a $K$ rational point, we can
assume that this element $Q$ is in the Jacobian and defined over
$L$.~\qed

Our goal here is to study elements of the Brauer group of the
product of the two curves $C$ and $C'$. To this end, we now add
the assumption that $L = F(C'')$ where $C'' \to C'$ is a cyclic
unramfied cover of degree $n$ and $C''$ has an $F$ rational point.
Note that $C'' \to C'$ induces a surjective homomophism $\Jac(\bar
C'') \to \Jac(\bar C')$, and the Galois group of this latter cover
is translation by elements in the cyclic kernel of order $n$. Let
$\sigma \in \Gal(L/K) = \Gal(F(C'')/F(C'))$ be a generator. Then
$\sigma: C'' \to C''$ induces a covering map $\sigma': \Jac(C'')
\to \Jac(C'')$ and $\sigma'$ is induced by addition of an order
$n$ element $P' \in \Jac(C'')$. Since $P' = \sigma'(0)$, $P'$ is $F$
rational.

The point $Q$ of Lemma 4.3 is $f': \Spec(F(C'')) \to \Spec(F(C''))
\times_F \Jac(C)$ and $\Spec(F(C'')) \to C''$ induces $$f':
\Spec(F(C'')) \to C'' \times_F \Jac(C).$$ Let $D$ be the closure
of the image of $f'$, and consider the induced map $D \to C''$
which birationally is the identity. Since $C''$ is the unique
non-singular model in $F(C'')$, and there is a desingularization
$D' \to D$, it follows that $D' = D$ and $D \to C''$ is the
identity. That is, $D$ is the graph of a morphism $$g:C'' \to
\Jac(C).$$ By the universal property of Jacobians, this induces a
homomorphism $$g:\Jac(C'') \to \Jac(C).$$

Since $Q$ is an $L$ point of $\Jac(C)$, it makes sense to form
$\sigma(Q)$ which is also a graph of a morphism and in fact
$\sigma(Q)$ is the graph of $g' = g \circ \sigma^{-1}$. From
$$\sigma(Q) - Q = P$$ we deduce that as morphisms $g'(x) = g(x) +
P$. But $\sigma$ on $\Jac(C'')$ is translation by~$P'$, so as
endomorphisms of $\Jac(C'')$ we have $g(x - P') = g(x) + P$, or
$$g(-P') = P.$$ We have shown:

\proclaim Proposition 4.4. In the above situation, $\alpha$ is split if and only if there is a homomorphism $g:\Jac(C'') \to \Jac(C)$ such that $g(P') = P$ where $P' \in \Jac(C'')$ generates the kernel of
$\Jac(C'') \to \Jac(C')$.

Let us note a consequence of the above in the case $C = E$ and
$C' = E'$ are both elliptic curves.

\proclaim Lemma 4.5. If $E$ and $E'$ are not isogenous, then
$\alpha$ as above is not split.

\proof Of course $C''$ is also an elliptic curve $E'',$ and $E''
\to E'$ is an isogeny. We can identify $E$, $E'$ and $E''$ with
their Jacobians and, of course, $E'$ and $E''$ are isogenous. If
$\alpha$ is split, the morphism $g: E'' \to E$ must be an
isogeny.~\qed

Note that Lemma~4.5 is a generalization of [C, p.~138] which states the
non-splitting in the case $n = 2$ (but is easily generalized to all $n$).

Recall that in Lemma 4.2 we wrote $\alpha$ as
$\Delta(L/K,\sigma,f)$, where $L = L(C)$ and $K =K(C) $. If $F$
contains a primitive $n$ root of 1, then $L = K(h^{1/n})$. Also,
we have assumed $L/K$ is everywhere unramified so the $C'$ divisor
of $h$ has the form $nE$ where $E$ is a divisor on $C'$. Of course
$E$ then defines an $n$ torsion point on~$\Jac(C')$. Moreover
$\alpha$ is the symbol algebra $(h,f)$ of degree $n$ over $F(C'
\times C)$.

For our example, we specialize to the case
where $F$ is algebraically closed, and $C$ and $C'$ are elliptic
curves which we can identify with their Jacobians. Let $n= p$ be a
prime. Consider two non-isogenous elliptic curves $E$ and~$E'$.
Let $P'$ be an element of order $p$ on $E$ and put $E_2 =E/<P'>.$
Let $Q \in E$ be an independent (of $P'$) element of order $p$,
let $E_3 = E/<Q>$, and let $P$ be the image of $P'$. Then $F(E) =
F(E_2)(a_2^{1/p})$ and $pP = (a_3)$ so by the above $(a_2,a_3)$ is
a split algebra over $F(E_2 \times E_3)$. Similarly we use $E'$ to
define $E_1$ and $E_4$ and a split algebra $(a_1,a_4)$ over $F(E_1
\times E_4)$. Set $V_1 = E_1 \times E_2$ and $V_2 = E_3 \times
E_4$. Now we work over $K = F(V_1 \times V_2) = F(E_1 \times E_2
\times E_3 \times E_4)$ and we view all the $a_i$ as also being
elements of $K$. Set $D_1 = (a_1,a_2)$, a degree $p$ symbol
algebra over $F(V_1)$, and $D_2 = (a_3,a_4)$, a degree $p$ symbol
algebra over $F(V_2)$.

\proclaim Theorem 4.6. Both $D_i$ are division algebras, but $D_1
\otimes_F D_2$ has a zero divisor.

\proof Set $D_i'  = D_i \otimes_{F(V_i)} K$. It suffices to show
that the $D_i$ are division algebras but $D_1' \otimes_K D_2'$ is
not a division algebra. Since $E_1$ and $E_2$ are not isogenous,
$D_1 = (a_1, a_2)$ is a division algebra over $F(V_1)$ by 4.5.
Similarly, $D_2$ is a division algebra. However, over $K$,
$(a_2,a_3)$ is split so $a_2$ is a norm from $K(a_3^{1/p})$.
Similarly, $a_4$ is a norm from $K(a_1^{1/p})$. Thus,
$K((a_2a_4)^{1/p})$ is a subfield of both $D_1' = (a_1,a_2)_K$ and
$D_2' = (a_3,a_4)_K$. Hence $D_1' \otimes_K D_2'$ is not a
division algebra.~\qed
\bigskip
\leftline{References}
\medskip
\noindent
[BLR] Bosch, S, Lutkebohmert, W, Raynaud, N.;
{\it N\'eron Models}, Springer Verlag Berlin/Heidelberg/New York 1990.
\medskip
\noindent
 [C] Colliot-Th\'el\`ene, J.-L., ``Exposant et Indice D'Algebres
Simples Non Ramifi\'ees'', L'Enseignement Math. {\bf 48} (2002),
p. 127--146.
\medskip
\noindent
[DI] DeMeyer, F and Ongraham, E.; {\it Separable Algebras over Commutative Rings},
Springer-Verlag Berlin/Heidelberg/New York (LNM \#181) 1971.
\medskip
\noindent
[Sc] Schacher, M. Problem 26 p.~380 in:
Gordon, R, (ed.) {\it Ring Theory, proceedings} Academic Press, New York 1971.
\medskip
\noindent
[Ho] Hoobler, R, ``A cohomological interpretation of the Brauer group of rings'',
Pacific Journal of Math {\bf 86} (1980), no. 1, 89--92.
\medskip
\noindent
[K] {\it Lectures on Resolutions of Singularities}, Princeton University Press,
Princeton/Oxford 2007 (Annals of Math Studies 166).
\medskip
\noindent
[M] Milne, {\it E\'tale Cohomology}; Princeton University Press Princeton, NJ 1980.
\medskip
\noindent
[P] Pirutka, A., ``A bound to kill the ramification over function fields'',
preprint, arXiv:1105.3942.
\medskip
\noindent
[Ro] Roquete, P., ``Splitting of Algebras by Function Fields in One Variable'',
Nagoya Math. J. {\bf 27} (1966) p. 625--642.
\medskip
\noindent
[S] ``Division Algebras over p-Adic Curves'' J. Ramanujan Math. Soc.,
{\bf 12} (1997), no. 1, p. 25-47 and ``Correction to Division Algebras over p-Adic
Curves'' J. Ramanujan Math. Soc. {\bf 13} (1998), no. 2, p.~125--129.
\medskip
\noindent
[S1]Saltman, D.J.,
``Division algebras over surfaces'', J. Algebra {\bf 320} (2008), no. 4, 1543-1585.
\medskip
\noindent
[Se] Serre, J.-P.; {\it Local Fields}, Springer-Verlag, New York/Heidelberg/Berlin 1979
\end